\RequirePackage{ifpdf}
\ifpdf 
\documentclass[pdftex]{sigma}
\else
\documentclass{sigma}
\fi

\begin{document}

\allowdisplaybreaks

\renewcommand{\thefootnote}{$\star$}

\renewcommand{\PaperNumber}{086}

\FirstPageHeading

\ShortArticleName{Natural Intrinsic Geometrical Symmetries}

\ArticleName{Natural Intrinsic Geometrical Symmetries\footnote{This paper is a
contribution to the Special Issue ``\'Elie Cartan and Dif\/ferential Geometry''. The
full collection is available at
\href{http://www.emis.de/journals/SIGMA/Cartan.html}{http://www.emis.de/journals/SIGMA/Cartan.html}}}

\Author{Stefan HAESEN~$^\dag$ and Leopold VERSTRAELEN~$^\ddag$}

\AuthorNameForHeading{S.~Haesen and L.~Verstraelen}

\Address{$^\dag$~Simon Stevin Institute for Geometry, Wilhelminaweg
1, 2042 NN Zandvoort, The Netherlands}
\EmailD{\href{mailto:stefan.haesen@geometryinstitute.org}{stefan.haesen@geometryinstitute.org}}

\Address{$^\ddag$~Katholieke Universiteit Leuven, Department of
Mathematics, Celestijnenlaan 200B bus 2400,\\
\hphantom{$^\ddag$}~B-3000 Leuven, Belgium} \EmailD{\href{mailto:leopold.verstraelen@wis.kuleuven.be}{leopold.verstraelen@wis.kuleuven.be}}

\ArticleDates{Received April 08, 2009, in f\/inal form August 25, 2009;  Published online September 02, 2009}

\Abstract{A proposal is made for what could well be \textit{the
most natural symmetrical Riemannian spaces which are homogeneous
but not isotropic}, i.e.\ of what could well be the most natural
class of symmetrical spaces beyond the spaces of constant
Riemannian curvature, that is, beyond the spaces which are
homogeneous and isotropic, or, still, the spaces which satisfy the
axiom of free mobility.}

\Keywords{parallel transport; holonomy; spaces of constant curvature; pseudo-symmetry}

\Classification{53A55; 53B20}

\section{Introduction}

Let $(M^{n},g)$ be an arbitrary Riemannian manifold, i.e.\ let
$M^{n}$ be any dif\/ferential manifold of dimension $n$ endowed with
any Riemannian metric $(0,2)$ tensor $g$. On the one hand, what
follows, with the usual care, can fully be developed also for
indef\/inite semi- or pseudo-Riemannian spaces, and, moreover, can
appropriately be extended to Cartan's ``generalized''-spaces. And,
on the other hand, when given some extra structure on $M^{n}$,
e.g.\ like a Kaehlerian or a Sasakian structure, what follows can
accordingly be well specif\/ied too. The present presentation
however will simply be restricted to some natural symmetries
occurring in intrinsic proper Riemannian geometry.

A geometrical symmetry of a Riemannian space $(M^{n},g)$ concerns
the invariance of some geometrical quantity of $(M^{n},g)$ under
the performance of some transformations def\/ined on $(M^{n},g)$
\cite{weyl}. Various types of symmetries in Riemannian geometry
can thus be considered, essentially depending on the kind of
quantities and of transformations in question. And, in complete
analogy with such \emph{intrinsic symmetries}, various types of
\emph{extrinsic symmetries} can be considered on submanifolds
$(M^{n},g)$ in Riemannian ambient spaces
$(\widetilde{M}^{n+m},\widetilde{g})$. By \emph{natural
geometrical symmetries}, we mean symmetries in the above sense for
which the quantities and the transformations involved are the most
natural indeed, at least in our opinion. As more technical as well
as more expository references for these intrinsic and for these
extrinsic symmetries, cf.~\cite{deszcz1,deszcz2,dillen1,haesen1,lumiste,verstraelen1}.


\section{Transformations on Riemannian manifolds}

The transformations on Riemannian manifolds $(M^{n},g)$ which will
be considered hereafter in our speculations about symmetry will be
\emph{the parallel transports fully around all the infinitesimal
co-ordinate parallelograms on these manifolds}. A 4-fold
motivation to devote some attention to these rather than maybe to
some other kinds of transformations is the following.

\medskip

$(i)$ Essentially, when studying \emph{symmetry}, one is interested
in \emph{the preservation of one or other sort of measures of some
geometrical ``beings''}, cf.\ the ``parlant'' of \'{E}lie Cartan,
\emph{after these beings have been transformed in one or other
manner}. Working in Riemannian geometry, these measures do most
naturally concern \emph{measures taken by the Riemannian metric
tensor}~$g$. And since this tensor may well change from point to
point on a general Riemannian manifold $(M^{n},g)$, it seems wise
to restrict to transformations on these manifolds which do bring
all beings living at any point $p\in M$ back to this same point
$p$, so that their measurements done before and after these
transformations are carried through are taken by the same measure
$g(p)$. This condition seems to be the most elementary one to
impose if one is willing to show respect for the metrical
structure of the space $(M^{n},g)$.

\medskip

$(ii)$ The Riemannian metrical structure $g$ is def\/ined on a
\emph{differential manifold} $M^{n}$. So, the transformations to
consider in the purpose of studying some basic symmetries of
Riemannian manifolds $(M^{n},g)$ should fundamentally involve the
essence of what is a dif\/ferential structure. An $nD$ dif\/ferential
structure on a space $M$ consists of an atlas of local co-ordinate
systems (or patches or  charts)
$x^{1},\ldots,x^{h},\ldots,x^{k},\ldots,x^{n}$ around its points
$p$, which, in case there exists between two of such systems a
non-trivial overlap, do have dif\/ferential transitions of their
co-ordinates there. So, the most elementary aspect of a
dif\/ferential structure on a space $M$ might well be the existence
of ``\emph{$x y$ co-ordinate parallelograms $\mathcal{P}$ cornered
at its points $p$''}, whereby we have renamed $x^{h}$ as $x$ and
$x^{k}$ as $y$ and, keeping all other co-ordinates around a point
$p=(x^{1},\ldots,x^{h-1},x,x^{h+1},\ldots,x^{k-1},y,x^{k+1},\ldots,x^{n})$
f\/ixed, consider the ``parallelogram'' $\mathcal{P}$ formed by the
co-ordinate lines $x^{h}=x$, $x^{k}=y$, $x^{h}=x+\Delta x$ and
$x^{k}=y+\Delta y$, for arbitrary increments~$\Delta x$ and~$\Delta y$ of the co-ordinates $x$ and $y$. Actually, this goes
back to the essence of desCartes' systematic use of planar
(oblique) $x y$ co-ordinates in his G\'{e}om\'{e}trie of 1637,
based there and then on the Euclidean parallel postulate.

\medskip

In the following, the natural basic tangent vectors
$\frac{\partial}{\partial x}(p) =\frac{\partial}{\partial
x^{h}}(p)$ and $\frac{\partial}{\partial
y}(p)=\frac{\partial}{\partial x^{k}}(p)$ at $p$ to the $x$- and
$y$-co-ordinate lines will be denoted by $\vec{x}$ and $\vec{y}$
respectively, and the tangent $2D$ plane, or, still, the $2D$
linear subspace of the tangent space $T_{p}M^{n}=\mathbb{R}^{n}$
to $M^{n}$ at $p$, which is spanned by $\vec{x}$ and $\vec{y}$
will be further denoted by $\overline{\pi}$:
$\overline{\pi}=\vec{x}\wedge\vec{y}$. And, combining $(i)$ and $(ii)$
it seems no more than appropriate when showing respect
respectively for the metrical and for the dif\/ferential structures
of general Riemannian manifolds $(M^{n},g)$, and, being in the
search of basic symmetries in Riemannian geometry, to consider
some kinds of transformations starting at arbitrary points $p$ of
$M$ and moving fully around co-ordinate parallelograms
$\mathcal{P}$ cornered at $p$ and tangent there to arbitrary
tangent 2-planes $\overline{\pi}=\vec{x}\wedge\vec{y}$.

\medskip

$(iii)$ Both in Riemann's and in Helmholtz's pioneering studies of
Riemannian spaces $(M^{n},g)$ \cite{helmholtz,riemann}, with their
inspirations coming (a.o.) from philosophical thoughts
primordially about physical space-time and about human vision,
respectively, the basic metrical geometrical structure, i.e.\
basically -in modern terminology- the metric tensor $g$, was
essentially def\/ined as an inf\/initesimal notion:
$g=g_{hk}\mbox{d}x^{h}\mbox{d}x^{k}$, $\big[$whereby the
components $g_{hk}$ of $g$, ($h,k\in\{1,2,\ldots,n\}$), are real
valued functions on (open parts of) the ``underlying''
dif\/ferential manifold, namely,
$g_{hk}:M\rightarrow\mathbb{R}\big]$. And this in itself
essentially was possible since working on a manifold $M^{n}$ to
which the inf\/initesimal calculus from $\mathbb{R}^{n}$ was
extended via mutually compatible local charts. Accordingly, after
thus ref\/lecting a bit deeper than in $(i)$ and $(ii)$ on the
inf\/initesimal characters of both the dif\/ferential structure and
the metrical structure of a Riemannian manifold, we will further
on consider transformations on $(M^{n},g)$ moving beings from any
point fully around all \emph{infinitesimal} co-ordinate
parallelograms $\mathcal{P}$ back to $p$, i.e., in the following,
$\Delta x$ and $\Delta y$ will be considered to be inf\/initesimal
quantities.

\medskip

$(iv)$ Finally, amongst all such transformations of beings around
such parallelograms $\mathcal{P}$, the most natural ones, again
with due respect for the dif\/ferential structure as well as for the
metrical structure of Riemannian manifolds $(M^{n},g)$, seem well
to be the corresponding \emph{parallel transports}. These were
introduced independently by Levi-Civita \cite{levicivita} and by
Schouten \cite{schouten}, in particular aiming for truly
geometrical insights in the meaning of the
\emph{Riemann--Christoffel curvature tensor} $R$ and of \emph{the
Riemann sectional curvatures} $K(p,\pi)$ of the $2D$ sections
$\pi$ of the tangent spaces $T_{p}M^{n}$ at the points $p$ of $M$.
These parallel transports of vectors along curves in $(M^{n},g)$
are conceptually equivalent with the \emph{Riemannian connection}
$\nabla$ which, in Koszul's approach, and as worked out by Nomizu,
according to the fundamental lemma of Riemannian geometry, is the
unique way to associate a vector f\/ield $\nabla_{X}Y$ with given
vector f\/ields $X$ and $Y$ on $(M^{n},g)$ following the derivation
rules of linearity and of Leibniz and being compatible with the
dif\/ferential structure and with the metrical structure of
$(M^{n},g)$, i.e., being \emph{symmetrical} (that is, the
commutator $\nabla_{X}Y-\nabla_{Y}X$ which is def\/ined for two such
vector f\/ields $X$ and $Y$ in terms of this connection coincides
with the commutator, or Lie bracket, $\big[X,Y\big] = X\, Y - Y\,
X$, which is def\/ined for these vector f\/ields by the dif\/ferential
structure) and being \emph{metrical} (that is, $Z\big[g(X,Y)\big]
- g(\nabla_{Z}X,Y)-g(X,\nabla_{Z}Y) = 0$, or, put otherwise,
$\nabla g=0$, or still, the scalar product of vectors does not
change under their parallel transports along curves).

\medskip

As is well known, and goes back to Schouten \cite{schouten}, one
has the following \emph{holonomy property}: after parallel
transport of a vector $\vec{z}\in T_{p}M^{n}$ fully around a
co-ordinate parallelogram $\mathcal{P}$ cornered at $p$ and
tangent there to $\overline{\pi}=\vec{x}\wedge\vec{y}$, results
the vector
\[ \vec{z}_{\overline{\pi}}^{\star} = \vec{z} + \big[
R(\vec{x},\vec{y})\vec{z}\big]\, \Delta x\Delta y +
\mathcal{O}^{>2}(\Delta x,\Delta y), \]
whereby $R(X,Y): TM\rightarrow TM$ is the
\emph{curvature operator}, i.e.\ $R(X,Y)Z :=
(\nabla_{X}\nabla_{Y}-\nabla_{Y}\nabla_{X}-\nabla_{[X,Y]})Z$, or
still, whereby $R: TM\times TM\times TM\rightarrow TM:
(X,Y,Z)\mapsto R(X,Y)Z$ is the $(1,3)$ Riemann--Christof\/fel
curvature tensor. Thus $R(\vec{x},\vec{y})\vec{z}$ measures the
second order change of a~vector $\vec{z}\in T_{p}M^{n}$ after
parallel transport around an inf\/initesimal co-ordinate
parallelogram~$\mathcal{P}$. Or, by taking into account that
parallel transports are isometries, and so, in particular, do not
change the lengths of vectors, one has the following.

\begin{proposition}[Schouten]
The curvature operators of Riemannian manifolds measure the
changes of directions at points under parallel transports fully
around the infinitesimal parallelograms cornered at these points.
\end{proposition}

Actually, the above holonomy property of Schouten during the last
decades is often used as the def\/inition of the $(1,3)$ curvature
tensor $R$. And, stated as the symmetry property of preserving
directions for the class of transformations under consideration,
one has the following.

\begin{theorem}[Schouten] \label{th:schouten}
The locally Euclidean $($or locally flat$)$ Riemannian manifolds $($the
spaces $(M^{n},g)$ for which $R\equiv 0)$, are precisely the
Riemannian manifolds for which all directions are invariant under
their parallel transports fully around all infinitesimal
co-ordinate parallelograms.
\end{theorem}

\begin{remark}
Actually, it was Cartan who very appropriately introduced the term
``holonomy'' by combining the Greek ``holos'' and ``nomos''.
\end{remark}

\section{Geometrical meaning of the metrical endomorphism}

The \emph{natural metrical endomorphism} $X\wedge_{g}Y:
TM\rightarrow TM$ associated with two vector f\/ields $X$ and~$Y$ on
a Riemannian manifold $(M^{n},g)$ is def\/ined by $(X\wedge_{g}Y)Z
:= g(Y,Z)X-g(X,Z)Y$. Let $\vec{x}$ and $\vec{y}$ be orthonormal
vectors at $p$, and let $\vec{z} = \vec{z}_{\overline{\pi}} +
\vec{z}_{\overline{\pi}^{\perp}}$ be the canonical orthogonal
decomposition of any vector~$\vec{z}$ at~$p$ in its components in
$\overline{\pi}=\vec{x}\wedge\vec{y}$ and in the orthogonal
complement $\overline{\pi}^{\perp}$ of $\overline{\pi}$ in
$T_{p}M^{n}=\mathbb{R}^{n}$. Now, we rotate
$\vec{z}_{\overline{\pi}}$ around $p$ in the plane
$\overline{\pi}$ over an inf\/initesimal angle $\Delta\varphi$, thus
obtaining a vector $(\vec{z}_{\overline{\pi}})_{\Delta\varphi}$,
and def\/ine the vector $\vec{z}_{\overline{\pi}}^{\wedge}
:=(\vec{z}_{\overline{\pi}})_{\Delta\varphi} +
\vec{z}_{\overline{\pi}^{\perp}}$. The procedure going from
$\vec{z}$ to $\vec{z}_{\overline{\pi}}^{\wedge}$ is called
\emph{the rotation of $\vec{z}$ at $p$ with respect to the plane
$\overline{\pi}$ over an angle $\Delta\varphi$}, and one has the
following:
\[ \vec{z}_{\overline{\pi}}^{\wedge} = \vec{z} +
\big[(\vec{x}\wedge_{g}\vec{y})\vec{z}\big]\, \Delta\varphi +
\mathcal{O}^{>1}(\Delta\varphi). \]
Thus the vector $(\vec{x}\wedge_{g}\vec{y})\vec{z}$
measures the f\/irst order change of the vector $\vec{z}$ after an
inf\/initesimal rotation of $\vec{z}$ at $p$ with respect to the
plane $\overline{\pi}=\vec{x}\wedge\vec{y}$, or, formulated more
loosely, we have the following.

\begin{proposition}[\cite{haesen1}]
The natural metrical endomorphisms $\wedge_{g}$ of Riemannian
manifolds measure the changes of directions at points under
infinitesimal rotations with respect to $2D$ planes at these points.
\end{proposition}

\section{The sectional curvatures of Riemann}

The $(0,4)$ \emph{Riemann--Christoffel curvature tensor} $R$ of
$(M^{n},g)$ is related to the $(1,3)$ tensor $R$ by
$R(X,Y,Z,W)=g(R(X,Y)Z,W)$. The simplest $(0,4)$ tensor which is
canonically determined on $(M^{n},g)$ and which has the same
algebraic symmetry properties as $R$ likely is the tensor
$G(X,Y,Z,W) := g((X\wedge_{g}Y)Z,W)$, or, still, in terms of the
Nomizu--Kulkarni product of $(0,2)$ tensors: $G =
\frac{1}{2}g\wedge g$, i.e., basically $G$ is the Nomizu--Kulkarni
square of the metric $(0,2)$ tensor~$g$.

Let $\vec{v}$ and $\vec{w}$ be any pair of linearly independent
vectors to $M$ at a point $p$ spanning there a $2D$ plane
$\pi=\vec{v}\wedge\vec{w}\subset T_{p}M^{n}$. Then, the
\emph{sectional curvature} $K(p,\pi) :=
R(\vec{v},\vec{w},\vec{w},\vec{v})/G(\vec{v},\vec{w},\vec{w},\vec{v})$
of~$M^{n}$ at $p$ for $\pi$, (the def\/inition is independent of the
choice of basis $\vec{v},\vec{w}$ for $\pi$), i.e., cf.\ Riemann,
the \emph{Gauss curvature} at $p$ of the 2-dimensional surface
$G_{\pi}^{2}\subset M^{n}$ consisting locally of \emph{all the
geodesics of $(M^{n},g)$ passing through $p$ and whose velocity
vectors at $p$ belong to $\pi$} (such that
$T_{p}G_{\pi}^{2}=\pi$), can be thought of as a kind of
nomalisation of the curvature operator by the natural metrical
endomorphism, or, still, of the parallel transport of directions
in $M$ at $p$ around co-ordinate parallelograms cornered at $p$ by
the rotations of directions in $M$ at $p$ with respect to planes
in $T_{p}M$. For any orthonormal basis $\vec{v}$, $\vec{w}$ of
$\pi$, of course: $K(p,\pi)=R(\vec{v},\vec{w},\vec{w},\vec{v})$.

Levi-Civita \cite{levicivita} gave a geometrical interpretation of
the sectional curvatures $K(p,\pi)$ in terms of his
\emph{parallelogramoids}, of which we next consider the special
case of his \emph{squaroids}. Let $\vec{v}$ and~$\vec{w}$ be any
orthonormal basis of a tangent 2-plane $\pi$ of $(M^{n},g)$ at~$p$. Consider the geodesic~$\alpha$ through $p=\alpha(0)$ with
velocity $\alpha'(0)=\vec{w}$, and on it then localise a point
$q=\alpha(\varepsilon)$ at an inf\/initesimal distance $\varepsilon$
from~$p$. Then move $\vec{v}$ parallel along $\alpha$ from $p$ to~$q$, thus obtaining at~$q$ the vector~$\vec{v}_{\alpha}^{\star}$.
Through~$p$ and~$q$, further, we consider the geodesics $\beta$
and $\gamma$, $\beta(0)=p$ and $\gamma(0)=q$, with respective
velocities $\beta'(0)=\vec{v}$ and
$\gamma'(0)=\vec{v}_{\alpha}^{\star}$. On these geodesics we then
localise the points $\overline{p}=\beta(\varepsilon)$ and
$\overline{q}=\gamma(\varepsilon)$ at distances $\varepsilon$ from
$p$ and $q$, respectively. Finally, a Levi-Civita squaroid $p\,
q\, \overline{q}\, \overline{p}$ is completed by joining
$\overline{p}$ and $\overline{q}$ by a geodesic $\delta$. In
general, the geodesic distance between~$\overline{p}$ and~$\overline{q}$ will be dif\/ferent from $\varepsilon$ and, denoting
this distance by $\varepsilon'$, one has the following.

\begin{theorem}[Levi-Civita]
In first order approximation: $K(p,\pi) = \big(\varepsilon^{2} -
\varepsilon^{\prime 2}\big)/\varepsilon^{4}$.
\end{theorem}

As references, and, in some cases, as sources for more precise
statements for the following classical results, see e.g.
\cite{kobayashi1,kobayashi2,kuhnel,wolf}.

\begin{theorem}[Cartan] \label{th:cartan}
The knowledge of the $(0,4)$ curvature tensor $R$ is equivalent to
the knowledge of all sectional curvatures $K(p,\pi)$ of Riemann.
\end{theorem}

Theorem~\ref{th:schouten} can be reformulated as follows.

\begin{theorem}[Schouten]
The Riemannian manifolds with constantly vanishing sectional
curvature of Riemann, $K(p,\pi)\equiv 0$, are precisely the spaces
for which all directions are invariant under their parallel
transports fully around all infinitesimal co-ordinate
parallelograms.
\end{theorem}

A Riemannian manifold $(M^{n},g)$ is said to be \emph{a space of
constant curvature} $K$ or \emph{a real space form} $M^{n}(c)$
when for all points $p$ and for all planar sections $\pi\subset
T_{p}M$ one has $K(p,\pi) = K = c \in\mathbb{R}$.

\begin{theorem}
A Riemannian manifold $(M^{n},g)$, $n\geq 2$, is a real space form
if and only if $R=c  G$, $c\in\mathbb{R}$.
\end{theorem}

\begin{theorem}[Schur]
A Riemannian manifold $(M^{n},g)$, $n\geq 3$, for which all
sectional curvatures~$K$ are isotropic, i.e., a Riemannian
manifold for which for every one of its points $p$ the curvatures
$K(p,\pi)$ are the same for all possible planar sections
$\pi\subset T_{p}M$, has constant sectional curvatures $K(p,\pi)$,
or, still, their sectional curvatures $K$ further are also
independent of the points $p$.
\end{theorem}

\begin{theorem}[Riemann]
Any two $nD$ Riemannian manifolds with the same constant sectional
curvatures are locally isometric; actually, in local co-ordinates,
the metrical fundamental form of a real space form of curvature
$c$ can be expressed as $\mathrm{d}s^{2} =
\left\{1+\frac{c}{4}\sum_{j}(x^{j})^{2}\right\}^{-2}
\sum_{i}(\mathrm{d}x^{i})^{2}$.
\end{theorem}

\begin{theorem}[Killing, Hopf]
The Euclidean spaces $\mathbb{E}^{n}$ ($c=0$) and the classical
non-Euclidean geometries, of elliptic type on the spheres
$\mathbb{S}^{n}$ $(c>0)$ and of hyperbolic type on the spaces
$\mathbb{H}^{n}$ $(c<0)$, are the model spaces of the real space
forms $M^{n}(c)$, $n\geq 2$.
\end{theorem}

Moreover, we recall that the real space forms may well be
considered to be the utmost possible symmetrical Riemannian
spaces, in the sense that they are \emph{homogeneous} (i.e.\ they
``behave'' the same at all of their points $p$) and that they are
\emph{isotropic} (i.e.\ they ``behave'' the same in all
directions). Both conditions together are a way to express that
the isometry groups of the real space forms are of maximal
dimension possible amongst all Riemannian manifolds of the same
dimensions, or, still, that we have the following.

\begin{theorem}[Riemann, Helmholtz, Lie, Klein, Tits]
A Riemannian manifold $(M^{n},g)$ is a real space form $M^{n}(c)$
if and only if it satisfies the axiom of free mobility.
\end{theorem}

Finally, we recall also the fundamental following result (for a
new recent proof of which we refer to Matveev \cite{matveev}),
which in some sense unif\/ies the Euclidean and the classical
non-Euclidean geometries from \emph{a projective point of view}.

\begin{theorem}[Beltrami] \label{th:beltrami}
The real space forms constitute the projective class of the
locally Euclidean spaces, or, still, by applying geodesic
transformations to locally Euclidean spaces one obtains spaces of
constant curvature and the class of the spaces of constant
curvature is closed under geodesic transformations.
\end{theorem}

\section{Geometrical meaning of semi-symmetry}

The curvatures $R$ and $K$ of a Riemannian manifold $(M^{n},g)$
are its \emph{main metrical invariants}~\cite{berger1,berger2}.
They essentially involve the \emph{second order derivatives}, as
those are determined by the dif\/ferential structure of the manifold
$M^{n}$, of the metrical structure of the Riemannian manifold,
i.e.~of~$g$.

Let $\vec{v}$ and $\vec{w}$ be any orthonormal vectors at any
point $p$ in a Riemannian manifold, and let~$\mathcal{P}$ be any
inf\/initesimal co-ordinate parallelogram cornered at $p$ as before,
i.e.\ in particular being tangent at $p$ to the plane
$\overline{\pi}=\vec{x}\wedge\vec{y}$. By parallel transport of
$\vec{v}$ and $\vec{w}$ around $\mathcal{P}$ then result the
orthonormal vectors $\vec{v}^{\star} =
\vec{v}+\big[R(\vec{x},\vec{y})\vec{v}\big]\, \Delta x\Delta y +
\mathcal{O}^{>2}(\Delta x,\Delta y)$ and $\vec{w}^{\star} =
\vec{w}+\big[R(\vec{x},\vec{y})\vec{w}\big]\, \Delta x\Delta y +
\mathcal{O}^{>2}(\Delta x,\Delta y)$, such that the plane
$\pi=\vec{v}\wedge\vec{w}$ at $p$ is accordingly parallel
transported around~$\mathcal{P}$ into the plane
$\pi^{\star}=\vec{v}^{\star}\wedge\vec{w}^{\star}$ at $p$. Hence,
$K(p,\pi^{\star}) =
R(\vec{v}^{\star},\vec{w}^{\star},\vec{w}^{\star},\vec{v}^{\star})
= R(\vec{v},\vec{w},\vec{w},\vec{v}) + \big[
R(R(\vec{x},\vec{y})\vec{v},\vec{w},\vec{w},\vec{v})
+R(\vec{v},R(\vec{x},\vec{y})\vec{w},\vec{w},\vec{v})
+R(\vec{v},\vec{w},R(\vec{x},\vec{y})\vec{w},\vec{v})
+R(\vec{v},\vec{w},\vec{w},R(\vec{x},\vec{y})\vec{v})\big]\,
\Delta x\Delta y + \mathcal{O}^{>2}(\Delta x,\Delta y) = K(p,\pi)
- \big[(R\cdot
R)(\vec{v},\vec{w},\vec{w},\vec{v};\vec{x},\vec{y})\big]\, \Delta
x\Delta y + \mathcal{O}^{>2}(\Delta x, \Delta y)$, whereby $R\cdot
R$ denotes the $(0,6)$ tensor which is obtained by the actions of
the curvature operators $R(X,Y)$ as derivations on the $(0,4)$
Riemann--Christof\/fel curvature tensor.

\begin{proposition}[\cite{haesen1}]
The curvature tensor $R\cdot R$ measures the changes of the
sectional curvatures $K(p,\pi)$ of a Riemannian manifold
$(M^{n},g)$ for all planes $\pi$ at all points $p$ under the
parallel transports of these planes $\pi$ fully around all
infinitesimal parallelograms $\mathcal{P}$ cornered at~$p$ and
tangent there to all planes~$\overline{\pi}$.
\end{proposition}

The Riemannian manifolds $(M^{n},g)$ for which $R\cdot
R=0$ are called \emph{semi-symmetric} or \emph{Szab\'{o}
symmetric}. These spaces were classif\/ied in general by Szab\'{o}
\cite{szabo1,szabo2}; for more specif\/ic information on some
particular cases, see also Boeckx \cite{boeckx}, Kowalski
\cite{kowalski1} and Lumiste \cite{lumiste}.

\begin{theorem}[\cite{haesen1}]
The semi-symmetric Riemannian spaces are precisely the Riemannian
manifolds for which, up to second order, all the sectional
curvatures are invariant under their parallel transports fully
around all infinitesimal co-ordinate parallelograms.
\end{theorem}

\begin{theorem}[Szab\'{o}] \label{th:szabo}
Let $(M^{n},g)$ be a semi-symmetric Riemannian manifold of
dimension $n>2$. Then there is an everywhere dense open subset
$U\subset M$ such that around every point $p\in U$ the space $M$
is locally isometric to a direct product of an open part of a
Euclidean space and some infinitesimally irreducible simple
semi-symmetric leaves $N$ which are: $(i)$ locally symmetric if
$\nu_{p}=0$ and $u_{p}>2$; $(ii)$ locally isometric to a Euclidean,
an elliptical or a hyperbolic cone if $\nu_{p}=1$ and $u_{p}>2$;
$(iii)$ locally isometric to a Kaehlerian cone if $\nu_{p}=2$ and
$u_{p}>2$; and $(iv)$ locally isometric to a space foliated by
Euclidean $2$-codimensional leaves if $\nu_{p}=n-2$ and $u_{p}=2$,
whereby $\nu_{p}$ and $u_{p}$ respectively are the indices of
nullity and of conullity of $(M^{n},g)$ at $p$.
\end{theorem}

The curvature tensor $R$ of a 2-dimensional Riemannian manifold
$(M^{2},g)$ essentially reduces to the Gauss curvature function
$K:M^{2}\rightarrow \mathbb{R}$ and so, by the Theorem of
Schwarz--Young, every $2D$ Riemannian space is automatically
semi-symmetric. The curvature condition $R\cdot R=0$ f\/irst
appeared in the studies concerning locally symmetric spaces by
P.A.~Shirokov and by \'{E}.~Cartan, namely as integrability
condition of the curvature condition $\nabla R=0$. For more
information about the origins of the notions of semi- and of
pseudo-symmetry, in particular concerning the mayor inf\/luences in
these contexts of the articles \cite{nomizu} and \cite{chen1} of
Nomizu and Chen, respectively, see \cite{lumiste} and
\cite{verstraelen1}. At this stage, we will conf\/ine to the comment
that whereas the curvature condition to be locally symmetric or
Cartan symmetric ($\nabla R\equiv 0$) essentially involve the 3rd
order derivatives of the metric of a Riemannian space $(M^{n},g)$,
the curvature condition to be semi-symmetric or Szab\'{o}
symmetric ($R\cdot R\equiv 0$) again essentially involves the 2nd
order derivatives of the metric $g$.

\begin{proposition}[\cite{haesen1}]
The $(0,6)$ tensor
\begin{gather*}
(R\cdot R)(X_{1},X_{2},X_{3},X_{4};X,Y)   :=   (R(X,Y)\cdot
R)(X_{1},X_{2},X_{3},X_{4})  \\
\qquad{}  =   -R(R(X,Y)X_{1},X_{2},X_{3},X_{4})
-R(X_{1},R(X,Y)X_{2},X_{3},X_{4}) \\
\qquad\quad{}  -R(X_{1},X_{2},R(X,Y)X_{3},X_{4})
-R(X_{1},X_{2},X_{3},R(X,Y)X_{4})
\end{gather*}
has the following algebraic symmetry properties:
\begin{gather*}
a) \ \ (R\cdot R)(X_{1},X_{2},X_{3},X_{4};X,Y) = -(R\cdot
R)(X_{2},X_{1},X_{3},X_{4};X,Y) \\
\qquad \quad {}=(R\cdot
R)(X_{3},X_{4},X_{1},X_{2};X,Y),
\\
b) \ \ (R\cdot R)(X_{1},X_{2},X_{3},X_{4};X,Y) +(R\cdot
R)(X_{1},X_{3},X_{4},X_{2};X,Y) \\
\qquad\quad{}+(R\cdot
R)(X_{1},X_{4},X_{2},X_{3};X,Y) =0,\\
c) \ \ (R\cdot R)(X_{1},X_{2},X_{3},X_{4};X,Y)= -(R\cdot
R)(X_{1},X_{2},X_{3},X_{4};Y,X),
\\
d) \ \ (R\cdot R)(X_{1},X_{2},X_{3},X_{4};X,Y) +(R\cdot
R)(X_{3},X_{4},X,Y;X_{1},X_{2}) \\
\qquad\quad {}+(R\cdot
R)(X,Y,X_{1},X_{2};X_{3},X_{4}) =0.
\end{gather*}
\end{proposition}


\section{The sectional curvatures of Deszcz}

The simplest non-trivial $(0,6)$ tensor which is canonically
determined on $(M^{n},g)$ and which has the same algebraic
symmetry properties as $R\cdot R$ likely is the Tachibana tensor
\begin{gather*}
(\wedge_{g}\cdot R)(X_{1},X_{2},X_{3},X_{4};X,Y)   :=
((X\wedge_{g}Y)\cdot R)(X_{1},X_{2},X_{3},X_{4}) \\
\qquad{}  =   -R((X\wedge_{g}Y)X_{1},X_{2},X_{3},X_{4})
-R(X_{1},(X\wedge_{g}Y)X_{2},X_{3},X_{4}) \\
\qquad\quad{} -R(X_{1},X_{2},(X\wedge_{g}Y)X_{3},X_{4})
-R(X_{1},X_{2},X_{3},(X\wedge_{g}Y)X_{4}),
\end{gather*}
i.e.\ is the $(0,6)$ tensor resulting from the actions
as derivations of the metrical endomorphisms $X\wedge_{g}Y$ on the
$(0,4)$ curvature tensor $R$.

\begin{theorem}[cf.~\cite{eisenhart}]
The Tachibana tensor of a Riemannian manifold $(M^{n},g)$, $n\geq
3$, vanishes identically, $\wedge_{g}\cdot R\equiv 0$, if and only
if $(M^{n},g)$ is a real space form $M^{n}(c)$.
\end{theorem}

In view of the geometrical interpretation of the
metrical endomorphism $\vec{x}\wedge_{g}\vec{y}$ given above in
terms of inf\/initesimal rotations at $p\in M^{n}$ with respect to
planes $\overline{\pi}=\vec{x}\wedge\vec{y}$, in particular, for
orthonormal vectors $\vec{v}$ and $\vec{w}$ at $p$ one has
$\vec{v}_{\overline{\pi}}^{\wedge} = \vec{v} +
\big[(\vec{x}\wedge_{g}\vec{y})\vec{v}\big]\, \Delta\varphi +
\mathcal{O}^{>1}(\Delta\varphi)$ and
$\vec{w}_{\overline{\pi}}^{\wedge} = \vec{w} +
\big[(\vec{x}\wedge_{g}\vec{y})\vec{w}\big]\, \Delta\varphi +
\mathcal{O}^{>1}(\Delta\varphi)$, such that the plane
$\pi=\vec{v}\wedge\vec{w}$ at $p$ is accordingly rotated at $p$
with respect to the plane $\overline{\pi}=\vec{x}\wedge\vec{y}$
into the plane
$\pi^{\wedge}=\vec{v}_{\overline{\pi}}^{\wedge}\wedge\vec{w}_{\overline{\pi}}^{\wedge}$
at $p$. And hence, completely analogously to a previous
calculation: $K(p,\pi^{\wedge}) = K(p,\pi) + \big[
(\wedge_{g}\cdot
R)(\vec{v},\vec{w},\vec{w},\vec{v};\vec{x},\vec{y})\big]\,
\Delta\varphi + \mathcal{O}^{>1}(\Delta \varphi)$.

\begin{proposition}[\cite{haesen1}]
The Tachibana tensor $\wedge_{g}\cdot R$ measures the changes of
the sectional curvatures $K(p,\pi)$ of a Riemannian manifold
$(M^{n},g)$ for all planes $\pi$ at all points $p$ under the
infinitesimal rotations of these planes $\pi$ at $p$ with respect
to all planes $\overline{\pi}$ at $p$.
\end{proposition}

\begin{theorem}[\cite{haesen1}]
The real space forms are precisely the Riemannian manifolds for
which, up to first order, all the sectional curvatures are
invariant under their infinitesimal rotations at all points with
respect to all planes.
\end{theorem}

Next, in analogy with the kind of normalisation of the changes of
directions under parallel transports around parallelograms
$\mathcal{P}$ cornered at~$p$ by the changes of directions under
rotations at~$p$ with respect to planes $\overline{\pi}$ at $p$,
which results in the notion of the \emph{sectional curvatures
$K(p,\pi)$ of Riemann}, we introduce the notion of \emph{the
double sectional curvatures} or \emph{the sectional curvatures
$L(p,\pi,\overline{\pi})$ of Deszcz} on Riemannian manifolds
$(M^{n},g)$ by now considering the changes of the sectional
curvatures $K(p,\pi)$ instead of the previous changes of
directions. In order to do so, for a trivial technical reason, we
can meaningly in this respect only consider planes
$\pi=\vec{v}\wedge\vec{w}$ and
$\overline{\pi}=\vec{x}\wedge\vec{y}$ at $p$ such that
$(\wedge_{g}\cdot
R)(\vec{v},\vec{w},\vec{w},\vec{v};\vec{x},\vec{y})\neq 0$, in
which case $\pi$ is said to be \emph{curvature dependent} on
$\overline{\pi}$. For such pairs of planes $\pi$ and
$\overline{\pi}$ then, their sectional curvature of Deszcz is
def\/ined by $L(p,\pi,\overline{\pi}) = (R\cdot
R)(\vec{v},\vec{w},\vec{w},\vec{v};\vec{x},\vec{y})/
(\wedge_{g}\cdot
R)(\vec{v},\vec{w},\vec{w},\vec{v};\vec{x},\vec{y})$. This
def\/inition is independent of the choices of bases for the planes
$\pi$ and $\overline{\pi}$. Similar to Theorem~\ref{th:cartan} of
Cartan, we have the following result and its consequences.

\begin{theorem}[\cite{haesen1}]
The knowledge of the $(0,6)$ curvature tensor $R\cdot R$ is
equivalent to the knowledge of all sectional curvatures
$L(p,\pi,\overline{\pi})$ of Deszcz.
\end{theorem}

\begin{theorem}[\cite{haesen1}]
A Riemannian manifold has constantly vanishing sectional
curvatures of Deszcz, $L(p,\pi,\overline{\pi})\equiv 0$, if and
only if it is semi-symmetric, i.e.\ if, up to second order, all its
sectional curvatures $K(p,\pi)$ of Riemann are invariant under the
parallel transport fully around all infinitesimal co-ordinate
parallelograms cornered at $p$ and tangent there to
$\overline{\pi}$.
\end{theorem}

In terms of the squaroids of Levi-Civita, a geometrical
interpretation of the sectional curvatures
$L(p,\pi,\overline{\pi})$ of Deszcz can be given as follows.
Consider a squaroid $\mathcal{S}$ of sides $\varepsilon$,
$\varepsilon'$, which is constructed upon orthonormal vectors
$\vec{v}$ and $\vec{w}$ at a point $p$. Then, let
$\vec{v}^{\star}$, $\vec{w}^{\star}$ and, respectively
$\vec{v}^{\wedge}_{\overline{\pi}}$,
$\vec{w}^{\wedge}_{\overline{\pi}}$ denote the orthonormal vectors
at $p$ which result, respectively, from the parallel translation
of $\vec{v}$, $\vec{w}$ around an inf\/initesimal co-ordinate
parallelogram $\mathcal{P}$ and from an associated inf\/initesimal
rotation of $\vec{v}$, $\vec{w}$ at $p$ with respect to a plane
$\overline{\pi}=\vec{x}\wedge\vec{y}$, whereby this association
means that the inf\/initesimal orders of the increments of the
co-ordinates or angles concerned do correspond, namely that
$\Delta\varphi=\Delta x\Delta y$. The sides of the squaroids
$\mathcal{S}^{\star}$ and $\mathcal{S}^{\wedge}$ constructed
respectively upon the vectors $\vec{v}^{\star}$, $\vec{w}^{\star}$
and $\vec{v}^{\wedge}_{\overline{\pi}}$,
$\vec{w}^{\wedge}_{\overline{\pi}}$  will be denoted respectively
by $\varepsilon$, $\varepsilon^{\star\prime}$ and by~$\varepsilon$,~$\varepsilon^{\wedge\prime}$.

\begin{theorem}[\cite{jahanara1}]
In first order approximation: $L(p,\pi,\overline{\pi}) =
\big[(\varepsilon^{\star\prime})^{2} -
(\varepsilon')^{2}\big]/\big[(\varepsilon^{\wedge\prime})^{2} -
(\varepsilon')^{2}\big]$.
\end{theorem}


\section{Geometrical meaning of pseudo-symmetry}

A Riemannian manifold $(M^{n},g)$, ($n\geq 3$), is said to be
\emph{pseudo-symmetric in the sense of Deszcz} or is called
\emph{Deszcz symmetric} if, for some function $L_{R}:M\rightarrow
\mathbb{R}$, $R\cdot R = L_{R}\, \wedge_{g}\cdot R$. We observe
that there does not hold a strict analog of the Theorem of Schur
(concerning the curvature function $K:M\rightarrow\mathbb{R}$ in
$R = K\, G$) in the case of pseudo-symmetric manifolds; see
\cite{deszcz1,deszcz3,haesen2} for examples for which the function
$L_{R}$ is \emph{not constant}. Following Kowalski and Sekizawa
\cite{kowalski2,kowalski3}, the pseudo-symmetric spaces with
constant function $L_{R}$ are said to be \emph{pseudo-symmetric of
constant type}.

\begin{theorem}[\cite{haesen1}]
A Riemannian manifold $(M^{n},g)$, $(n\geq 3)$, is Deszcz
symmetric if and only if its sectional curvature of Deszcz
$L(p,\pi,\overline{\pi})$ is isotropic, i.e., if
$L(p,\pi,\overline{\pi})$ is independent of the planes $\pi$ and
$\overline{\pi}$, or, still, if the double sectional curvature
function $L(p,\pi,\overline{\pi})$ actually is a~function
$L=L_{R}:M\rightarrow \mathbb{R}$.
\end{theorem}

The condition for a Riemannian manifold to be
\emph{pseudo-symmetric}, a terminology which f\/irst appeared, as
far as we know, in \cite{deszczgrycak}, did occur in Grycak's
investigations of semi-symmetric warped products \cite{grycak}, as
well as in the study of geodesic mappings on semi-symmetric
spaces, amongst others by Sinyukov, Mike\v{s} and Venzi
\cite{mikesh,venzi}. From these studies and from the artic\-le~\cite{defever} of Defever and Deszcz we quote the following
theorem, to be considered in some sense in analogy with
Theorem~\ref{th:beltrami} of Beltrami.

\begin{theorem}[Sinyukov, Mike\v{s}, Venzi, Defever and Deszcz]
If a semi-symmetric Riemannian space admits a geodesic
transformation onto some other Riemannian manifold, then this
latter manifold must itself be pseudo-symmetric, and, if a
pseudo-symmetric Riemannian space admits a geodesic transformation
onto some other Riemannian manifold, then this latter manifold
must itself also be pseudo-symmetric.
\end{theorem}

In any case, in the late 19seventies and early 19eighties the
relevance of the intrinsic pseudo-symmetry became more clear
mainly by some studies concerning the geometry of submanifolds,
notably starting with the studies, in particular by Deszcz, on
\emph{extrinsic spheres in semi-symmetric spaces}, which extended
studies on totally umbilical submanifolds in Cartan symmetric
spaces by Chen \cite{chen1} and Olszak \cite{olszak}.


\section{Further types of pseudo-symmetry curvature conditions}

In full analogy with the above mentioned studies concerning the
parallel transport around inf\/initesimal co-ordinate parallelograms
of the Riemann curvatures leading to the notions of
pseudo-symmetry of Riemannian manifolds and of the sectional
curvatures of Deszcz, the consi\-de\-ra\-tion of the Weyl conformal
curvatures and of the Ricci curvatures instead lead to the notions
of \emph{Weyl and Ricci pseudo-symmetry} and of the \emph{Weyl and
Ricci curvatures of Deszcz}~\mbox{\cite{jahanara1,jahanara2}}. For some
insight into the relationships between these various types of
pseudo-symmetry curvature conditions, see e.g.~\cite{deszcz2,marrakesh}.

The following seems to be a further natural curvature condition of
a similar kind, $C\cdot C = L_{C}\, \wedge_{g}\cdot C$, whereby
$C$ denotes the $(0,4)$ Weyl conformal curvature tensor as well as
the corresponding curvature operator and $L_{C}$ is a real
function on the manifold $M^{n}$, which condition, in contrast to
the above conditions, is \emph{invariant under conformal
transformations} of the metric $g$. Spaces $(M^{n},g)$ satisfying
this condition are said to have \emph{a pseudo-symmetric Weyl
tensor C}.

In connection with further structures on Riemannian spaces, the
above mentioned studies can be adapted accordingly. For instance,
on a Kaehlerian manifold $(M^{n},g, J)$ it seems most natural to
focus in particular on the invariance of the holomorphic sectional
curvatures under their parallel transport around inf\/initesimal
holomorphic co-ordinate parallelograms, etc.; the results on
several of such specialisations of the above general Riemannian
pseudo-symmetry curvature conditions are being considered at
present.


\section{Some examples of pseudo-symmetric spaces}

Of course, according to Theorem~\ref{th:szabo}, starting from
Szab\'{o}'s classif\/ication of the semi-symmetric Riemannian
spaces, by applying, eventually iteratively, geodesic
transformations, one always obtains Riemannian manifolds
$(M^{n},g)$ which are pseudo-symmetric in the sense of Deszcz.

On the other hand, let $M^{n}$ be a hypersurface in a Euclidean
space $\mathbb{E}^{n+1}$, $n\geq 3$. Amongst the simplest possible
forms of the shape operator of $M^{n}$ in $\mathbb{E}^{n+1}$, one
has those whereby the principal curvatures at every point are (1):
$(0,0,\ldots,0)$; (2): $(\lambda,\lambda,\ldots,\lambda)$,
$\lambda\neq 0$; (3): $(\lambda,0,\ldots,,0)$, $\lambda\neq 0$;
(4):~$(\lambda,\ldots,\lambda,0,\ldots,0)$, $\lambda\neq 0$ and
$\lambda$ appearing more than once; (5):
$(\lambda,\mu,0,\ldots,0)$, $\lambda\neq 0\neq \mu$ and
$\lambda\neq \mu$; (6): $(\lambda,\mu,\ldots,\mu)$, $\lambda\neq
0\neq\mu$ and $\lambda\neq\mu$; (7):
$(\lambda,\ldots,\lambda,\mu,\ldots,\mu)$, $\lambda\neq 0\neq\mu$
and $\lambda\neq\mu$ and both $\lambda$ and $\mu$ appearing more
than once. Then there are the following correspondences:
$M^{n}\subset\mathbb{E}^{n+1}$ is totally geodesic in case (1);
(non-totally geodesic) totally umbilical in case~(2); (non-totally
geodesic) cylindrical in case~(3); the cases (1), (2) and~(3)
together cover the hypersurfaces of constant sectional curvature,
i.e.\ the $M^{n}$ in $\mathbb{E}^{n+1}$ which are real space forms
(which for these hypersurfaces is equivalent to being Einstein)
and cases (1) and (3) deal with the locally f\/lat hypersurfaces;
semi-symmetric hypersurfaces which are not real space forms
concern the cases (4) and (5), so that the semi-symmetric
hypersurfaces~$M^{n}$ of $\mathbb{E}^{n+1}$ correspond to (1),
(2), (3), (4) and (5), as shown by Nomizu; conformally f\/lat
$M^{n}$, for $n> 3$, which are not of constant curvature
correspond to case (6); and the intrinsic pseudo-symmetric~$M^{n}$
in $\mathbb{E}^{n+1}$ correspond to (1), (2), (3), (4), (5), (6)
and (7). So, \emph{a hypersurface $M^{n}$ in $\mathbb{E}^{n+1}$ is
a non semi-symmetric, intrinsically pseudo-symmetric Riemannian
manifold if and only if it has exactly two non-zero principal
curvatures $\lambda$ and $\mu$, and then its double sectional
curvature is given by $L = \lambda\, \mu$}. A few further notes
might be in place here. $(i)$~A~hypersurface $M^{n}$ in
$\widetilde{M}^{n+1}$ is said to be \emph{quasi-umbilical} when it
has a principal curvature of multiplicity $\geq n-1$, and,
\emph{in conformally flat ambient spaces}, as shown by Cartan and
Schouten, \emph{this is equivalent to $M^{n}$ being itself
conformally flat}, whenever $n>3$. $(ii)$ We recall that a
hypersurface $M^{n}$ in $\widetilde{M}^{n+1}(c)$ is called a
\emph{cyclide of Dupin} if it has, at every point, exactly two
distinct principal curvatures which are both constant along the
corresponding principal tangent subbundles, see e.g.~\cite{cecil1,cecil2}. $(iii)$ A hypersurface $M^{n}$ in
$\widetilde{M}^{n+1}$ is said to be \emph{$2$-quasi-umbilical} when
it has a~principal curvature of multiplicity $\geq n-2$, and, in
particular, a hypersurface $M^{n}$ in Euclidean space
$\mathbb{E}^{n+1}$, with dimension $n>3$, \emph{has a
pseudo-symmetric Weyl tensor} if and only if $M$ is
\emph{pseudo-symmetric} or \emph{its principal curvatures are
given by $(\lambda,\mu,\nu,\ldots,\nu)$ whereby
$\lambda\neq\mu\neq\nu\neq\lambda$}.

A submanifold $M^{n}$ of a Riemannian manifold
$\widetilde{M}^{n+m}$ is said to be \emph{$\vec{H}$-parallel} if
\mbox{$R^{\perp}\cdot\vec{H} = \vec{0}$}, i.e.\ if, up to second order
approximation, the mean curvature vector $\vec{H}$ remains
invariant under the $\nabla^{\perp}$ parallel transport along $M$
in $\widetilde{M}$ completely around inf\/initesimal co-ordinate
parallelograms~$\mathcal{P}$ in~$M$~\cite{dillen1}. \emph{All
pseudo-parallel submanifolds $M$ in $\widetilde{M}$ are
$\vec{H}$-parallel}, and, \emph{every extrinsically
pseudo-symmetric submanifold $M^{n}$ in a real space form
$\widetilde{M}^{n+m}(\widetilde{c})$ is either minimal or
pseudo-umbilical or has flat normal connection, and if the first
normal spaces of $M^{n}$ in $\widetilde{M}^{n+m}$ have maximal
dimension $n(n+1)/2$ then $M$ is minimal or pseudo-umbilical in
$\widetilde{M}$}. All parallel, semi-parallel, and respectively,
pseudo-parallel submanifolds $M^{n}$ in real space forms
$\widetilde{M}^{n+m}(\widetilde{c})$ are automatically
intrinsically locally symmetric, semi-symmetric, and,
respectively, pseudo-symmetric Riemannian manifolds. For
classif\/ication results on these submanifolds, a.o.\ due to Ferus,
Dillen and N\"{o}lker, Lumiste, Asperti, Lobos, Tojeiro, Mercuri,
i.p.\ see \cite{asperti,dillennolker,ferus,lobos,lumiste}.

Recently, \emph{the conjecture on the Wintgen inequality}, as
formulated in \cite{desmet} was completely resolved in the
af\/f\/irmative by Lu~\cite{lu} and by Ge and Tang~\cite{ge}.

\begin{theorem}[Lu, Ge and Tang]
The Wintgen inequality $\rho\leq H^{2} -\rho^{\perp}
+\widetilde{c}$ holds for every submanifold $M^{n}$ in any real
space form $\widetilde{M}^{n+m}(\widetilde{c})$, whereby $\rho$
and $\rho^{\perp}$ are the normalised scalar curvature of $M$ and
the normalised scalar normal curvature of $M$ in $\widetilde{M}$,
respectively, and $H^{2}$ is the squared mean curvature of $M$ in
$\widetilde{M}$, and $\rho = H^{2}-\rho^{\perp} + \widetilde{c}$
if and only if, with respect to suitable tangent and normal
orthonormal frames $\{E_{i}\}$ and $\{\xi_{\alpha}\}$, the shape
operators of $M^{n}$ in $\widetilde{M}^{n+m}(\widetilde{c})$
assume the forms
\begin{gather*}
A_{1} = \left[\begin{array}{@{\,\,}c@{\,\,}c@{\,\,}c@{\,\,}c@{\,\,}c@{\,\,}}
              \lambda +\mu\cos\theta &  0 &  0 & \cdots & 0 \\
              0 & \lambda-\mu\cos\theta & 0 & \cdots & 0 \\
              0 & 0 & \lambda & \cdots & 0 \\
              \vdots & \vdots & \vdots & \ddots & \vdots \\
              0 & 0 & 0 & \cdots & \lambda \\
              \end{array} \right], \qquad
A_{2} = \left[\begin{array}{@{\,\,}c@{\,\,}c@{\,\,}c@{\,\,}c@{\,\,}c@{\,\,}}
              \mu\sin\theta & 0 & 0 & \cdots & 0 \\
              0 & -\mu \sin\theta & 0 & \cdots & 0 \\
              0 & 0 & 0 & \cdots & 0 \\
              \vdots & \vdots & \vdots & \ddots & \vdots \\
              0 & 0 & 0 & \cdots & 0 \\
              \end{array} \right],
\\ A_{3} = \left[\begin{array}{@{\,\,}c@{\,\,}c@{\,\,}c@{\,\,}c@{\,\,}c@{\,\,}}
              0 & \mu & 0 & \cdots & 0 \\
              \mu & 0 & 0 & \cdots & 0 \\
              0 & 0 & 0 & \cdots & 0 \\
              \vdots & \vdots & \vdots & \ddots & \vdots \\
              0 & 0 & 0 & \cdots & 0 \\
              \end{array} \right],\qquad A_{4}=\cdots=A_{m}=0.
\end{gather*}

{\noindent}Hereby $\rho =
\frac{2}{n(n-1)}\sum_{i<j}R(E_{i},E_{j},E_{j},E_{i})$, and
$\rho^{\perp} = \frac{2}{n(n-1)}\left[
\sum_{i<j}\sum_{\alpha<\beta}R^{\perp}(E_{i},E_{j};\xi_{\alpha},\xi_{\beta})^{2}\right]^{1/2}$.
\end{theorem}

There are many submanifolds satisfying the equality
$\rho=H^{2}-\rho^{\perp}+\widetilde{c}$, and, in analogy with
Chen's nomenclature for the submanifolds satisfying the equality
in several other kinds of general inequalities between intrinsic
and extrinsic invariants of submanifolds, these submanifolds are
called \emph{Wintgen ideal submanifolds}. For explicit
descriptions of Wintgen ideal submanifolds and for further
discussions on this inequality, one can consult amongst others
studies by Boruvka~\cite{boruvka}, Bryant
\cite{rysz:Bryant,bryant}, Choi--Lu \cite{choi},
Dajczer--Tojeiro \cite{dajczer1,dajczer2},
Dillen--Fastenakels--Van der Veken \cite{dillen4,dillen5}, Eisenhart~\cite{eisenhart2}, Guadalupe--Rodriguez \cite{guadalupe}, Kommerell
\cite{kommerell}, Rouxel \cite{rouxel} and Wintgen~\cite{wintgen}.
As was shown in \cite{deszcz99,petrovic1,petrovic2}, there are
close connections between Wintgen ideal submanifolds and intrinsic
pseudo-symmetry conditions. For instance, \emph{for $n>3$, every
Wintgen ideal submanifold~$M^{n}$ in
$\widetilde{M}^{n+m}(\widetilde{c})$ has pseudo-symmetric
conformal Weyl tensor}, and the \emph{minimal Wintgen ideal
submanifolds are characterised by the fact that} $L_{C} =
-\frac{n-3}{(n-1)(n-2)} K_{\inf}$; moreover, \emph{the
Deszcz symmetric Wintgen ideal submanifolds are either totally
umbilical}, (in particular, then $M$ being a real space form),
\emph{or minimal} (in which case $M$ is pseudo-symmetric of
constant type, namely $L=\widetilde{c}$).

Finally, we will comment on the fact that the pseudo-symmetric
Riemannian manifolds $(M^{n},g)$, $n\geq 3$, could well be
considered as being the most natural symmetric spaces which, in
extension of the perfect symmetry behaviour of the real space
forms, do admit privileged directions at all of its points, or,
still, as being \emph{the most symmetric anisotropic spaces}. In
particular, for spaces of dimension $n=3$, one has the following.

\begin{theorem}[\cite{deszczver}]\label{th:3D}
A $3D$ Riemannian manifold is Deszcz symmetric if and only if, it
is quasi-Einstein.
\end{theorem}

We recall that $(M^{n},g)$ is said to be an
\emph{Einstein} manifold when its Ricci tensor $S$ is proportional
to the metric $g$, or, still, when all its Ricci curvatures at all
points are equal, and $(M^{n},g)$ is said to be
\emph{quasi-Einstein} when at all points it has a Ricci curvature
of multiplicity $\geq n-1$; in particular, $(M^{n},g)$ is properly
quasi-Einstein when at all points it has a Ricci curvature of
multiplicity precisely $n-1$ and the other Ricci curvature then of
course has multiplicity~1. And, in this situation, in particular,
the Ricci principal direction with respect to this latter Ricci
curvature then determines on $(M^{n},g)$ a tangent direction which
geometrically is essentially dif\/ferent from all other tangent
directions. We recall the following.

\begin{theorem}[Schouten and Struik]
A $3D$ Riemannian manifold is an Einstein space if and only if it is
a real space form.
\end{theorem}

The real space forms $M^{n}(c)$ are semi-symmetric and
thus in particular pseudo-symmetric. So the proper $3D$ Deszcz
symmetric spaces are the proper quasi-Einstein $3D$ spaces, with a~Ricci curvature $\lambda$ of multiplicity 1 and a~Ricci curvature
$\mu\neq \lambda$ of multiplicity 2, and
$L=\lambda/2: $ $M\rightarrow\mathbb{R}$. In his approach to
geometries as topological manifolds $M$ endowed with
transformation groups~$G$ satisfying by def\/inition just a few
conditions such as to allow for the class of real space forms to
be enlarged so as to moreover incorporate certain still
homogeneous, but essentially anisotropic spaces, Thurston
\cite{thurston1,thurston2} introduced his so-called \emph{model
geometries} which can be seen as to determine successfully some of
the most natural anisotropic geometrical manifolds $(M^{n},G)$,
$n\geq 3$, beyond the projective class of the Euclidean and the
classical non-Euclidean geometries. The introduction of the
pseudo-symmetric spaces, i.e.\ of the projective class of the
semi-symmetric spaces, is a metrical approach also to determine
some of the most natural anisotropic manifolds through the
symmetry property given by the invariance of the curvature of the
connection $\nabla$ on a dif\/ferential manifold $M$ under the
parallel transport corresponding to $\nabla$ around the
inf\/initesimal co-ordinate parallelograms on $M$. We recall that
the eight $3D$ Thurston model geometries are $\mathbb{E}^{3}$,
$\mathbb{S}^{3}$, $\mathbb{H}^{3}$,
$\mathbb{S}^{2}\times\mathbb{E}^{1}$,
$\mathbb{H}^{2}\times\mathbb{E}^{1}$, the Heisenberg group
$H_{3}$, the covering group of the special linear group
$\widetilde{SL}(2,\mathbb{R})$ and the solvable Lie group $Sol$. A
canonical metric can be placed on each of the model spaces, which,
except for $\mathbb{H}^{3}$ and $Sol$, can be written as
$\mbox{d}s^{2} = \frac{d x^{2}+d y^{2}}{[1+m(x^{2}+y^{2})]^{2}} +
\left[d z + \frac{l}{2}\frac{y\, d x-x\, d
y}{1+m(x^{2}+y^{2})}\right]^{2}$, whereby if $m=l=0$,
$M=\mathbb{E}^{3}$; $4m-l^{2}=0$, $M=\mathbb{S}^{3}$; if $m<0$ and
$l=0$, $M=\mathbb{H}^{2}\times\mathbb{E}^{1}$; if $m>0$ and $l=0$,
$M=\mathbb{S}^{2}\times\mathbb{E}^{1}$; if $m<0$ and $l\neq 0$,
$M=\widetilde{SL}(2,\mathbb{R})$; and if $m=0$ and $l\neq 0$,
$M=H_{3}$. On the Lie group $Sol$, one can put the metric
$\mbox{d}s^{2}=\mbox{e}^{2z}\mbox{d}x^{2}+\mbox{e}^{-2z}\mbox{d}y^{2}
+\mbox{d}z^{2}$. Besides the real space forms $\mathbb{E}^{3}$,
$\mathbb{S}^{3}$ and $\mathbb{H}^{3}$, the f\/ive other geometries
are proper quasi-Einstein, and so one has the following.

\begin{theorem}[\cite{belkhelfa}]
All $3D$ Thurston geometries are either spaces of constant curvature
or Deszcz symmetric of constant type; $($for $\mathbb{E}^{3}$,
$\mathbb{S}^{3}$ and $\mathbb{H}^{3}$, $K=0$, $1$ and $-1$
respectively; for $\mathbb{S}^{2}\times\mathbb{E}^{1}$ and
$\mathbb{H}^{2}\times\mathbb{E}^{1}$, $L=0$; for $H_{3}$ and
$\widetilde{SL}(2,\mathbb{R})$, $L=1$; and for $\mbox{Sol}$,
$L=-1)$.
\end{theorem}

Deleting $\mbox{Sol}$ from the $3D$ Thurston spaces and
adding instead $SU(2)$, with the above metric whereby $m>0$ and
$l\neq 0$, one obtains the list of the $3D$ so-called \emph{d'Atri
spaces}, i.e.\ the Riemannian manifolds for which locally all
geodesic ref\/lections are volume preserving. One has the following.

\begin{theorem}[\cite{belkhelfa}]
All $3D$ d'Atri spaces are Deszcz symmetric of constant type
$($besides the above, for $SU(2)$, $L=1)$.
\end{theorem}

Basically, Theorem~\ref{th:3D} readily follows from the fact that
for every Riemannian manifold of dimension 3 the Weyl conformal
curvature tensor $C$ automatically vanishes. In the above
respects, for higher dimensions one has the following.

\begin{theorem}[\cite{deprez,deszcz4}]
A conformally flat Riemannian manifold $(M^{n},g)$ of dimension
$n\geq 4$ is Deszcz symmetric if and only if it has at most two
distinct Ricci curvatures $($of arbitrary multiplicities$)$.
\end{theorem}

For the Einstein spaces $S=\lambda\, g$ this again
actually only concerns the real space forms, and in case the
spaces under consideration are not Einstein, the two orthogonally
complementary eigenspaces of their Ricci tensor consist of
geometrically non-equivalent tangent directions, thus realising on
these pseudo-symmetric spaces a manifest, yet in some sense
elementary, anisotropy.

As was already mentioned in the Introduction, the above
essentially goes through for semi-Riemannian spaces. In
particular, in \cite{haesen2} a classif\/ication, based on the
algebraic properties of the Weyl and Ricci tensor, was obtained
for the 4-dimensional pseudo-symmetric space-times. It follows
that most of the well-known space-times, such as the Schwarzschild
metric, the Reissner--N\"{o}rdstrom metric, the Kottler metric and
the Friedmann--Lema\^{\i}tre--Robertson--Walker metrics, are
pseudo-symmetric but not semi-symmetric. Moreover, the f\/irst three
metrics are examples of non-conformally f\/lat, pseudo-symmetric
spaces. See also \cite{defever1,deszcz6} for more details.


\subsection*{Acknowledgements}
The authors do thank the referees whose comments resulted in real
improvements of the original version of this paper.  The f\/irst author was partially supported by the Spanish MEC Grant
MTM2007-60731 with FEDER funds and the Junta de Andaluc\'{\i}a
Regional Grant P06-FQM-01951. Both authors were partially
supported by the Research Foundation Flanders project G.0432.07.


\pdfbookmark[1]{References}{ref}
\LastPageEnding

\end{document}